\documentclass[fleqn]{article}
\usepackage[utf8]{inputenc}
\usepackage{enumerate,amsthm,amssymb,xcolor,hyperref,comment}

\usepackage[leqno]{amsmath}

\makeatletter
\newcommand{\leqnomode}{\tagsleft@true}
\newcommand{\reqnomode}{\tagsleft@false}
\makeatother

\usepackage{latexsym}
\usepackage{amsfonts}
\usepackage[leqno]{amsmath}
\usepackage{tikz}
\usetikzlibrary{decorations.pathmorphing}
\tikzset{snake it/.style={decorate, decoration=snake}}
\usepackage{float}
\usepackage{lmodern}

\def\dd{\hbox{-}}

\usepackage{marvosym}               
\def\longbox#1{\parbox{0.85\textwidth}{#1}}

\begin{document}

\title{A note on simplicial cliques}
\author{Maria Chudnovsky\thanks{Supported by NSF grant DMS-1763817.}\\
Princeton University, Princeton, NJ 08544, USA
\\
\\
Alex Scott\\
Mathematical Institute, University of Oxford, Oxford OX2 6GG, UK
\\
\\
Paul Seymour\thanks{Supported by AFOSR grant A9550-19-1-0187 and NSF grant DMS-1800053.}\\
Princeton University, Princeton, NJ 08544, USA
\\
\\
Sophie Spirkl\thanks{We acknowledge the support of the Natural Sciences and Engineering Research
  Council of Canada (NSERC), [funding reference number RGPIN-2020-03912].
  Cette
recherche a été financée par le Conseil de recherches en sciences naturelles
 et en génie du Canada (CRSNG), [numéro de référence RGPIN-2020-03912]
   .}\\
University of Waterloo, Waterloo, ON N2L 3E9, Canada}

\date {November 5, 2020; revised \today}

\newtheorem{theorem}{}[section]

\maketitle
\begin{abstract}
Motivated by an application in condensed matter physics and  quantum information theory, we prove that every non-null even-hole-free claw-free graph has a simplicial clique, that is, a clique $K$ such that for every
  vertex $v \in K$, the set of neighbours of $v$ outside of $K$ is a clique.
In fact, we prove the existence of a simplicial clique in a more general class of graphs defined by forbidden induced subgraphs.
\end{abstract}

\section{Introduction}

All graphs in this paper are finite and simple.
Let $G$ be a graph. A {\em hole} in $G$ is an
induced cycle of length at least four. By a {\em path} we always mean an
induced path. For $v \in V(G)$ we denote by $N(v)$
the set of neighbours of $v$ (so $v \not \in N(v)$).
$G$ is {\em even-hole-free} if all holes in $G$ have odd length,
and $G$ is {\em claw-free} if $G$ has no induced subgraph isomorphic to
$K_{1,3}$.   A non-empty set $K \subseteq V(G)$ is a {\em simplicial clique} if
$K$ is a clique, and for every $v \in K$ we have that $N(v) \setminus K$
is a clique. The unique element of a  simplicial clique  of size one is called a {\em simplicial vertex}.

This paper is motivated by a question from condensed matter physics and
quantum information theory concerning so-called {\em spin models}, i.e. models of interacting qubits (two-level quantum systems).
Each spin model is defined by a Hamiltonian operator, and to every such
Hamiltonian
one can associate a graph, called its {\em frustration graph}. In
\cite{models} a new method is given that allows us  to ``solve a model''
(meaning in this case to find the spectrum and the eigenvectors of
the Hamiltonian) whose frustration graph is even-hole-free,
claw-free, and has a simplicial clique. This augments earlier results
of \cite{linegraphs} where it is shown that models whose frustration
graphs are line-graphs are solvable using certain classical tools. 
The solution method of \cite{models} uses only the structure of the
frustration graph, and it is an extension of both \cite{Fendley} and
\cite{Lieb}.  The authors of \cite{models} raised a question:

\begin{theorem} \label{question1}
  {\bf Question:} Does every non-null even-hole-free claw-free graph have a simplicial clique?
\end{theorem}

  In other words, does their new solvability result hold for all models whose
frustration graphs are even-hole-free and claw-free? In this note we answer
their question affirmatively (the ``dome of an edge'' is defined before
the statement of \eqref{subgraphs}); in fact we prove a stronger result:

\begin{theorem} \label{simplclique}
  Let $G$ be a non-null even-hole-free claw-free graph.
    \begin{enumerate}
  \item If $G$ is chordal, then  $G$ has a simplicial vertex. 
  \item  For every hole $H$ of $G$ there is  an edge $ab$ of $H$ such that
    the dome of $ab$ is a simplicial clique.
 \end{enumerate}
   In particular,      $G$ has a simplicial  clique.
\end{theorem}

We have an even stronger
result, describing explicitly the structure of all such graphs, but the proof is harder, and we do not present it here. The main result of this paper is
\ref{subgraphs} which is a strengthening of \ref{simplclique}, and we will
explain it in Section \ref{stronger}.

We  remark that the answer to
\ref{question1} becomes negative if we omit either the assumption that the graph is even-hole-free or that the graph is  claw-free. The complement of a cycle of odd length at least seven is  a claw-free graph with no
simplicial clique.  Moreover, the square of a cycle of length at least nine is
an example of a $C_4$-free claw-free graph with no simplicial clique.
(The {\em square} of a graph $G$ is the graph obtained from $G$ by
making every vertex adjacent to all its second neighbours.)
And here is an example of an even-hole-free graph rather than just $C_4$-free.
Let $k$ be an odd positive integer.
The following is a construction of an  even-hole-free graph $G_k$ with $2k$ vertices and with no simplicial clique. Let  the vertex set of
$G_k$ be the union of $k$ disjoint  pairs of adjacent vertices $\{a_i,b_i\}$ where $i \in \{1, \ldots, k\}$. For $i=\{1, \ldots, k-1\}$ 
add the edges $a_ia_{i+1}, a_ib_{i+1},b_ia_{i+1}$; add also the edges
$a_ka_1,a_kb_1,b_ka_1$. There are no more edges in $G$. It is easy to check that
$G_k$ is even-hole-free and has no simplicial clique.

In \cite{oldalg} an algorithm is presented
that finds a simplicial clique in a claw-free graph if one exists.
The authors of \cite{models} also asked if that algorithm can be
simplified when the input is known to be  even-hole-free. An easy corollary
of our main result \ref{simplclique}  is such a simpler, but slower,
algorithm \ref{alg}. In fact \ref{alg} works under the more general assumptions
of \ref{subgraphs}.

\section{A strengthening} \label{stronger}

The goal of this section is to present our main result \ref{subgraphs}.

Let $G$ be a graph. For $X \subseteq V(G)$ we denote by $G[X]$ the
graph induced by $G$ on $X$.
For $A \subseteq V(G)$ and $x \in V(G) \setminus A$, we say that $x$
is {\em complete} to $A$ if $x$ is adjacent to every element of $A$, and that
$x$ is {\em anticomplete} to $A$ if $x$ is non-adjacent to every element of $A$.
Two disjoint subsets $A,B \subseteq V(G)$ are {\em complete } to each other if
every vertex of $B$ is complete to $A$, and {\em anticomplete} to each other
if every vertex of $B$ is anticomplete to $A$.

\begin{figure}[h]
\begin{center}
\begin{tikzpicture}[scale=0.28]

\node[inner sep=2.5pt, fill=black, circle] at (-2, 2)(v1){}; 
\node[inner sep=2.5pt, fill=black, circle] at (5, 2)(v2){}; 
\node[inner sep=2.5pt, fill=black, circle] at (0, 4)(v3){}; 
\node[inner sep=2.5pt, fill=black, circle] at (3, 4)(v4){};
\node[inner sep=2.5pt, fill=black, circle] at (1.5, 1)(v5){}; 

\draw[black, dotted, thick] (v1)  .. controls +(0,-5.5) and +(0,-5.5) .. (v2);
\draw[black, thick] (v1) -- (v3);
\draw[black, thick] (v2) -- (v4);
\draw[black, thick] (v3) -- (v4);
\draw[black, thick] (v3) -- (v5);
\draw[black, thick] (v1) -- (v5);
\draw[black, thick] (v2) -- (v5);
\draw[black, thick] (v4) -- (v5);


\node[inner sep=2.5pt, fill=black, circle] at (9, 4)(v1){}; 
\node[inner sep=2.5pt, fill=black, circle] at (12, 4)(v2){}; 
\node[inner sep=2.5pt, fill=black, circle] at (12, -2)(v3){}; 
\node[inner sep=2.5pt, fill=black, circle] at (9, -2)(v4){};
\node[inner sep=2.5pt, fill=black, circle] at (10.5, 1)(v5){}; 

\draw[black, dotted, thick] (v1)  .. controls +(-3,0) and +(-3,0) .. (v4);
\draw[black, dotted, thick] (v2)  .. controls +(3,0) and +(3,0) .. (v3);
\draw[black, thick] (v1) -- (v2);
\draw[black, thick] (v3) -- (v4);
\draw[black, thick] (v1) -- (v5);
\draw[black, thick] (v2) -- (v5);
\draw[black, thick] (v3) -- (v5);
\draw[black, thick] (v4) -- (v5);


\node[inner sep=2.5pt, fill=black, circle] at (16, 4)(v3){}; 
\node[inner sep=2.5pt, fill=black, circle] at (18, 1)(v4){};
\node[inner sep=2.5pt, fill=black, circle] at (16, -2)(v5){}; 
\node[inner sep=2.5pt, fill=black, circle] at (24, 4)(v6){}; 
\node[inner sep=2.5pt, fill=black, circle] at (22, 1)(v7){}; 
\node[inner sep=2.5pt, fill=black, circle] at (24, -2)(v8){};

\draw[black, thick] (v3) -- (v5);
\draw[black, thick] (v4) -- (v5);
\draw[black, thick] (v3) -- (v4);
\draw[black, thick] (v3) -- (v6);
\draw[black, dotted, thick] (v4) -- (v7);
\draw[black, thick] (v5) -- (v8);
\draw[black, thick] (v6) -- (v7);
\draw[black, thick] (v7) -- (v8);
\draw[black, thick] (v6) -- (v8);


\def\r{3.1}
\def\x{29}
\def\y{.9}
\node[inner sep=2.5pt, fill=black, circle] at ({\x+\r*cos(90)}, {\y+\r*sin(90)})(v1){}; 
\node[inner sep=2.5pt, fill=black, circle] at ({\x+\r*cos(90+360/7)}, {\y+\r*sin(90+360/7)})(v2){}; 
\node[inner sep=2.5pt, fill=black, circle] at ({\x+\r*cos(90+2*360/7)}, {\y+\r*sin(90+2*360/7)})(v3){};
\node[inner sep=2.5pt, fill=black, circle] at ({\x+\r*cos(90+3*360/7)}, {\y+\r*sin(90+3*360/7)})(v4){}; 
\node[inner sep=2.5pt, fill=black, circle] at ({\x+\r*cos(90+4*360/7)}, {\y+\r*sin(90+4*360/7)})(v5){}; 
\node[inner sep=2.5pt, fill=black, circle] at ({\x+\r*cos(90+5*360/7)}, {\y+\r*sin(90+5*360/7)})(v6){};
\node[inner sep=2.5pt, fill=black, circle] at ({\x+\r*cos(90-360/7)}, {\y+\r*sin(90-360/7)})(v7){};

\draw[black, thick] (v1) -- (v2);
\draw[black, thick] (v2) -- (v3);
\draw[black, thick] (v3) -- (v4);
\draw[black, thick] (v4) -- (v5);
\draw[black, thick] (v5) -- (v6);
\draw[black, thick] (v6) -- (v7);
\draw[black, thick] (v7) -- (v1);

\draw[black, thick] (v1) -- (v3);
\draw[black, thick] (v2) -- (v4);
\draw[black, thick] (v3) -- (v5);
\draw[black, thick] (v4) -- (v6);
\draw[black, thick] (v5) -- (v7);
\draw[black, thick] (v6) -- (v1);
\draw[black, thick] (v7) -- (v2);

\end{tikzpicture}
\end{center}
\vspace{-0.8cm}
\caption{A jewel, a line wheel, a short prism and a seven-antihole (here  dotted lines represent paths).}
\label{fig:clean}
\end{figure}
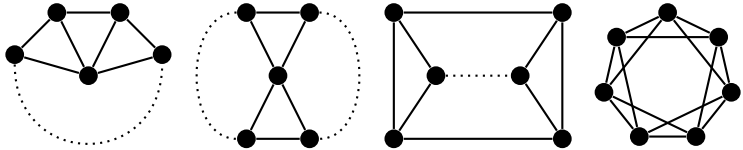

Next we define  a few types of graphs.
A graph is called {\em chordal}
if it has no holes.
A {\em jewel} is a graph consisting of a hole $H=h_1 \dd \ldots \dd h_k \dd h_1$
with $k \geq 4$ and a vertex $v \not \in V(H)$ such that
$N(v) \cap V(H)=\{h_1,h_2,h_3,h_4\}$.
A {\em line wheel} is a graph consisting of a hole
$H=h_1 \dd \ldots \dd h_k \dd h_1$ with $k \geq 6$ and a vertex $v \not \in V(H)$ such that there exists
$i \in \{4, \ldots, k-2\}$ with $N(v) \cap V(H)=\{h_1,h_2,h_i,h_{i+1}\}$.
A {\em short prism} is a graph consisting of a hole $h_1 \dd h_2 \dd h_3 \dd h_4 \dd h_1$ and a path $p_1 \dd \ldots \dd pk$ such that
$\{p_1, \ldots, p_k\} \cap \{h_1, h_2, h_3,h_4\}= \emptyset$,
$p_1$ is adjacent to $h_1$ and to $h_2$, and $p_k$ is adjacent to $h_3$ and to $h_4$,
and there are  no other edges between $\{p_1, \ldots, p_k\}$ and
$ \{h_1, h_2, h_3,h_4\}$.
Finally, the {\em seven-antihole}  is the complement of a cycle of seven
vertices. These graphs are depicted in Figure \ref{fig:clean}.

In what follows, whenever graph containment is mentioned, we
will mean containment as an induced subgraph.
We say that a graph $G$ is {\em clean} if  $G$ is claw-free and
contains no jewel, line wheel, short prism or seven-antihole. Note that clean
graphs may contain even holes.

First we show:
\begin{theorem}
  \label{clean}
  If $G$ is claw-free and even-hole-free, then $G$ is  clean.
  \end{theorem}

\begin{proof}
  Since $G$ is even-hole free, and in  particular $C_4$-free, $G$ does not contain short prisms or seven-antiholes.  Since a jewel contains a hole
  $H$ of length $k$, and a hole $h_1 \dd v \dd h_4 \dd h_5 \dd \ldots \dd h_k \dd h_1$ of length $k-1$, $G$ does not contain a jewel. Finally,
  at least one of the holes $H$, $h_2 \dd \ldots \dd h_i \dd v \dd h_2$,
  and $h_{i+1} \dd \ldots \dd h_1 \dd v \dd h_{i+1}$ is even, and so $G$
  does not contain a line wheel. This proves \ref{clean}.
  \end{proof}

We need one more definition. Let $ab$ be an edge of a graph $G$. Let
$X(ab)= \{a,b\} \cup(N(a) \cap N(b))$. The {\em dome} of $ab$ is the
set of all vertices $y \in X(ab)$ such that  $N(y) \setminus X(ab)$ is 
a clique. We call the set $X(ab) \setminus Y(ab)$ the {\em dome} of $ab$.
We can now state our main result.

\begin{theorem} \label{subgraphs}
  Let $G$ be a non-null clean graph.
  \begin{enumerate}
  \item If $G$ is chordal then $G$ has  a simplicial vertex.
  \item  For every  hole $H$ of $G$ there is  an edge
    $ab$ of $H$ such that
    the dome of $ab$ is a simplicial clique.
 \end{enumerate}
   In particular,      $G$ has a simplicial  clique.
\end{theorem}

In view of \ref{clean} we immediately deduce \ref{simplclique}.

\section{The proof of the main theorem}

In this section we prove \ref{subgraphs}.
We start with a lemma. 
\begin{theorem}
  \label{holenbrs}
  Let $G$ be a clean graph, let $H$ be a hole in $G$, and let $v \in V(G) \setminus V(H)$. Then one of the following holds:
  \begin{enumerate}
    \item $v$ is anticomplete to $V(H)$.
    \item $|V(H)|=5$ and $v$ is complete to $V(H)$.
    \item $v$ has exactly two neighbours in $H$ and they are consecutive.
      \item $v$ has exactly three neighbours in $H$ and they form a path of $H$.
  \end{enumerate}
\end{theorem}

\begin{figure}[h]
\begin{center}
\begin{tikzpicture}[scale=0.28]

  
\node[inner sep=2.5pt, fill=black, circle] at (-13, 2)(v1){}; 
\node[inner sep=2.5pt, fill=black, circle] at (-6, 2)(v2){}; 
\node[inner sep=2.5pt, fill=black, circle] at (-11, 4)(v3){}; 
\node[inner sep=2.5pt, fill=black, circle] at (-8, 4)(v4){};
\node[inner sep=2.5pt, fill=black, circle] at (-9.5,1)(v5){}; 

\node at (-9.5, 0) {$v$};

\draw[black, dotted, thick] (v1)  .. controls +(0,-5.5) and +(0,-5.5) .. (v2);
\draw[black, thick] (v1) -- (v3);
\draw[black, thick] (v2) -- (v4);
\draw[black, thick] (v3) -- (v4);

\def\r{3}
\def\x{0.5}
\def\y{1}

\node[inner sep=2.5pt, fill=black, circle] at ({\x+\r*cos(-90)}, {\y+\r*sin(-90)})(v1){}; 
\node[inner sep=2.5pt, fill=black, circle] at ({\x+\r*cos(-18)}, {\y+\r*sin(-18)})(v2){}; 
\node[inner sep=2.5pt, fill=black, circle] at ({\x+\r*cos(54)}, {\y+\r*sin(54)})(v3){};
\node[inner sep=2.5pt, fill=black, circle] at ({\x+\r*cos(126)}, {\y+\r*sin(126)})(v4){}; 
\node[inner sep=2.5pt, fill=black, circle] at ({\x+\r*cos(198)}, {\y+\r*sin(198)})(v5){}; 
\node[inner sep=2.5pt, fill=black, circle] at ({\x}, {\y})(v6){};

\node at (1, 0) {$v$};

\draw[black, thick] (v1) -- (v2);
\draw[black, thick] (v2) -- (v3);
\draw[black, thick] (v3) -- (v4);
\draw[black, thick] (v4) -- (v5);
\draw[black, thick] (v5) -- (v1);

\draw[black, thick] (v1) -- (v6);
\draw[black, thick] (v2) -- (v6);
\draw[black, thick] (v3) -- (v6);
\draw[black, thick] (v4) -- (v6);
\draw[black, thick] (v5) -- (v6);

\node[inner sep=2.5pt, fill=black, circle] at (7, 2)(v1){}; 
\node[inner sep=2.5pt, fill=black, circle] at (14, 2)(v2){}; 
\node[inner sep=2.5pt, fill=black, circle] at (9, 4)(v3){}; 
\node[inner sep=2.5pt, fill=black, circle] at (12, 4)(v4){};
\node[inner sep=2.5pt, fill=black, circle] at (10.5, 1)(v5){}; 

\node at (10.5,0) {$v$};

\draw[black, dotted, thick] (v1)  .. controls +(0,-5.5) and +(0,-5.5) .. (v2);
\draw[black, thick] (v1) -- (v3);
\draw[black, thick] (v2) -- (v4);
\draw[black, thick] (v3) -- (v4);
\draw[black, thick] (v3) -- (v5);
\draw[black, thick] (v4) -- (v5);

\node[inner sep=2.5pt, fill=black, circle] at (16, 2)(v1){}; 
\node[inner sep=2.5pt, fill=black, circle] at (23, 2)(v2){}; 
\node[inner sep=2.5pt, fill=black, circle] at (18, 4)(v3){}; 
\node[inner sep=2.5pt, fill=black, circle] at (21, 4)(v4){};
\node[inner sep=2.5pt, fill=black, circle] at (19.5, 1)(v5){}; 

\node at (19.5,0) {$v$};

\draw[black, dotted, thick] (v1)  .. controls +(0,-5.5) and +(0,-5.5) .. (v2);
\draw[black, thick] (v1) -- (v3);
\draw[black, thick] (v2) -- (v4);
\draw[black, thick] (v3) -- (v4);
\draw[black, thick] (v3) -- (v5);
\draw[black, thick] (v4) -- (v5);
\draw[black, thick] (v1) -- (v5);

\end{tikzpicture}
\end{center}
\vspace{-0.8cm}
\caption{Outcomes of \ref{holenbrs} (here  dotted lines represent paths).}
\label{fig:holenbrs}
\end{figure}
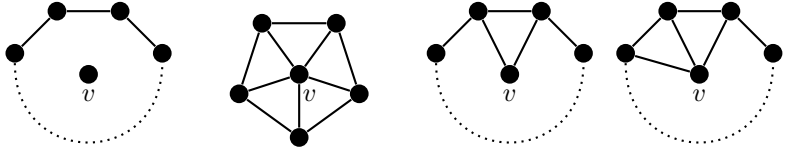

\begin{proof}
  The outcomes of \ref{holenbrs} are depicted in Figure \ref{fig:holenbrs}.
Write $H=h_1 \dd \ldots h_k \dd h_1$.
We may assume that $v$ has a neighbour in $V(H)$, for otherwise \ref{holenbrs}.1 holds. If $v$ is complete to $V(H)$,
then, since $G$ is claw-free and since  $G[V(H) \cup \{v\}]$ is not a jewel,
it follows that $k=5$ and so \ref{holenbrs}.2 holds.
  Thus we may assume that $v$ has a non-neighbour in $V(H)$, say $v$ is
  adjacent to $h_1$ and not to $h_k$. Since $G$ is claw-free, $v$ is adjacent to
  $h_2$, and $|N(v) \cap V(H)| \leq 4$.
  We may assume that $v$ has a neighbour in
  $V(H) \setminus \{h_1,h_2,h_3\}$, for otherwise \ref{holenbrs}.3 or
  \ref{holenbrs}.4 holds. Since $G$ is claw-free,
  $N(v) \cap (V(H) \setminus \{h_1,h_2,h_3\})$ is a clique, and therefore
  $|N(v) \cap (V(H) \setminus \{h_1,h_2,h_3\})| \leq 2$.
  Let $i \in \{4, \ldots,k-1 \}$ be minimum such that
  $v$ is adjacent to $h_i$. Since $\{v,h_{i-1}, h_i,h_{i+1}\}$ is not a
  claw, it follows that either
  \begin{itemize}
  \item $i=4$ and $v$ is adjacent to $h_3$, or
\item     $v$ is adjacent to $h_{i+1}$ and $v$ has no other neighbours in $V(H) \setminus \{h_1,h_2,h_i, h_{i+1}\}$.
  \end{itemize}
  In the former case $G[V(H) \cup \{v\}]$ is a jewel.
  Thus we may assume that the latter case holds. Since
  $|N(v) \cap V(H)| \leq 4$, it follows that $N(v) \cap V(H)=\{h_1,h_2,h_i,h_{i+1}\}$. But now   $G[V(H) \cup \{v\}]$ is
  a line wheel, a contradiction. This proves \ref{holenbrs}.
  \end{proof}

Now we turn to the proof of \ref{subgraphs}.

\begin{proof}
  If $G$ has no hole, then $G$ is chordal, and therefore has a simplicial vertex \cite{chordal}, and \ref{subgraphs} holds. Thus we may assume that $G$ has a
  hole.
    For an integer $k \geq 4$ a subset
  $W \subseteq V(G)$ is {\em $k$-structured} if  (here the addition is  mod $k$):
  \begin{itemize}
\item $W$ is the disjoint union of    $k$ non-empty cliques
  $K_1, \ldots, K_k$, 
\item for every $i \in \{1, \ldots, k\}$ every $v \in K_i$ has a neighbour in $K_{i-1}$ and a neighbour in $K_{i+1}$, and
\item if $i,j \in \{1, \ldots, k\}$ and $i \neq j \pm 1$
  then $K_i$ is anticomplete to $K_j$.
\end{itemize}
We call the partition $(K_1, \ldots, K_k)$ a {\em $k$-structure of} $W$.

Let $H$ be a hole of $G$. Then $H$ has length $k \geq 4$, and 
$V(H)$ is a $k$-structured set. If possible, choose $H$ with
$k \geq 5$.
Let $W \subseteq V(G)$ be a $k$-structured set with $k$-structure
$(K_1, \ldots, K_k)$, where each $K_i$ contains exactly one vertex of $H$,
and  such that $W$ is inclusion-wise maximal with this property.
In what follows addition and subtraction of indices of the $k$-structure is
mod $k$.

\vspace*{-0.4cm}
  
  \begin{equation} \label{opposite}
  \longbox{\emph{Let $i \in \{1, \ldots, k\}$. If $a,b \in K_i$ and $N(b) \cap K_{i+1} \not \subseteq N(a) \cap K_{i+1}$,  then $N(b) \cap K_{i-1} \subseteq N(a) \cap K_{i-1}$}.}
  \end{equation}
  We may assume $i=1$.
  If for each $j \in \{2,k\}$ there exists $a_j \in (N(b) \setminus N(a)) \cap K_j$,
  then $\{b,a_k,a,a_2\}$ is a claw in $G$, a contradiction. This proves
  \eqref{opposite}.

\vspace*{-0.4cm}

\begin{equation} \label{comparable}
  \longbox{\emph{Let $i \in \{1, \ldots, k\}$. For every $a,b \in K_i$ either
      $N(a) \cap K_{i+1} \subseteq N(b) \cap K_{i+1}$,  or $N(b) \cap K_{i+1} \subseteq N(a) \cap K_{i+1}$}.}
\end{equation}

We may assume $i=1$.
  Suppose there exist $a' \in (N(a) \setminus N(b)) \cap K_2$ and
  $b' \in (N(b) \setminus N(a)) \cap K_2$.  Since $K_1$ and $K_2$
  are cliques, $a$ is adjacent to $b$, and $a'$ is adjacent to $b'$.
  Now $a \dd a' \dd b' \dd b \dd a$ is a hole of length four.
  Let $C=N(a) \cap K_k$  and $C'= N(a') \cap K_3$.
  By \eqref{opposite} $b$  is complete to $C$, and $b'$ is complete
  to $C'$. Switching the roles of $a$ and $b$, we deduce that $b$ is
  anticomplete to $K_k \setminus C$, and, similarly, $b'$ is
  anticomplete to $K_3 \setminus C'$. Since the graph $G[\bigcup_{j=3}^kK_j]$
  is connected, there is a path $P=p_1 \dd \ldots \dd p_t$ in
  $G[\bigcup_{j=3}^kK_j]$ with $p_1 \in C$ and $p_t \in C'$; we may assume $P$ is
  chosen with $t$ minimum. Then $p_2, \ldots, p_{t-1} \not \in C \cup C'$,
  and therefore $V(P) \setminus \{p_1, p_t\}$ is anticomplete to
  $\{a,b,a',b'\}$. But now $G[a,a',b,b',p_1, \ldots, p_t]$ is a short
    prism in $G$, a contradiction.  This proves \eqref{comparable}.

\vspace*{-0.4cm}
\begin{equation} \label{complete}
  \longbox{\emph{For every $i \in \{1, \ldots, k\}$, $K_i$ is complete to at least one of $K_{i-1},K_{i+1}$}.}
\end{equation}

We may assume $i=1$.
  Suppose there exist $u \in K_k$, $v,w \in K_1$ and $z \in K_2$
  (where possibly $v=w$) such that $u$ is not adjacent to $w$,
  and $v$ is not adjacent to $z$.

  First we claim that $v,w$ can be chosen in such a way that $uv$ and
  $wz$ are edges. Suppose not; then we may assume that $v$ is non-adjacent
  to both $u$ and $z$. Since $(K_1, \ldots, K_k)$ is a $k$-structure, there
  exist $n_u,n_z \in K_1$ such that $u$ is adjacent to $n_u$, and
  $z$ is adjacent to $n_z$ (possibly $n_z=w$). Since $v$ is anticomplete to $\{u,z\}$,
  it follows from \eqref{opposite} that $n_u$ is non-adjacent to $z$,
  and $n_z$ is non-adjacent to $u$. But now we can choose $v=n_u$
  and $w=n_z$, and the claim holds.

  In view of the claim in the previous paragraph we assume that
  $uv$ and  $wz$ are adjacent (and in particular $v \neq w$).
  Let $u' \in K_k$ be a neighbour of $w$; by \eqref{comparable}
  $u'$ is adjacent to $v$. Now the set
  $T=\bigcup_{i=3}^{k-1}K_i$
  is non-empty  and connected, and  both $u'$ and $z$ have neighbours in $T$.
  Consequently, there is a path $R$ from $u'$ to $z$ with
  $V(R) \setminus \{u',z\} \subseteq T$. Let $x \in V(R)$ be the
  neighbour of $u'$. Then $x \in K_{k-1}$, and by \eqref{opposite} $x$ is
  adjacent to $u$.    It follows from the definition
  of a $k$-structure that $V(R) \setminus \{z,x,u'\}$ is anticomplete to
  $\{u,u',v,w\}$.
  But now the hole $z \dd w \dd v \dd u \dd x \dd R \dd z$ together with
  the vertex $u'$ forms a jewel in $G$, a contradiction.
  This proves
    \eqref{complete}.

\vspace*{-0.4cm}

\begin{equation} \label{path}
  \longbox{\emph{Let $j \in \{1, \ldots, k\}$. For every $i \in \{2, \ldots, k-2\}$,
      $a_j \in K_j$ and $a_{i+j} \in K_{j+i}$,
      there is a path $P$ from $a_j$ to $a_{j+i}$ with
      $V(P) \subseteq \bigcup_{t=j}^{j+i}K_{t}$ and using exactly one vertex from
      each of $K_j, \ldots, K_{j+i}$}.}
\end{equation}

We may assume $j=1$.
  The proof is by induction on $i$. Suppose first that $i=2$.
  Since by \eqref{complete} $K_2$ is complete to at least one of
  $K_1,K_3$, it follows that $a_1$ and $a_3$ have a common neighbour
  $a_2 \in K_2$. Now $a_1 \dd a_2 \dd a_3$ is the required path.

  Now assume that  $i>2$, let $a_{j+i-1} \in K_{j+i-1}$ be a neighbour of $a_{j+i}$.
  By the inductive hypothesis there is a path $P$ from $a_1$ to $a_{j+i-1}$
  with $V(P) \subseteq \bigcup_{r=j}^{i+j-1}K_r$ using exactly one vertex
  from      each of $K_j, \ldots, K_{i+j-1}$.  Now $a_j \dd P  \dd a_{j+i-1} \dd a_{j+i}$ is the
  required path. This proves \eqref{path}.



 \vspace*{-0.4cm}
\begin{equation} \label{claims}
  \longbox{\emph{Let $v \in V(G) \setminus W$. For $i \in \{1, \ldots, k\}$
      let $N_i=K_i \cap N(v)$ and  $M_i=K_i \setminus N_i$. The following hold
      for every $i \in \{1, \ldots, k\}$:
      \begin{enumerate}
      \item $N_i$ is anticomplete to at least one of $M_{i-1},M_{i+1}$.
        \item If $k>4$, then $M_i$ is anticomplete to at  least one of $N_{i-1},N_{i+1}$.
\end{enumerate}
    }}
\end{equation}

We may assume $i=1$. By \eqref{complete} we may assume that $K_1$ is complete to
$K_2$.

We first prove the first statement. We may assume that
there exists $m_2 \in M_2$, for otherwise the claim holds
($N_1$ is anticomplete to $M_2$ because $M_2=\emptyset$).
Now if $n_1 \in N_1$ has a neighbour $m_k \in M_k$, then
$\{n_1,v,m_2,m_k\}$ is a claw, a contradiction. This proves that
$N_1$ is anticomplete to $M_2$, and \eqref{claims}.1
follows.

Next we prove the second statement. We may assume that
there exist $m_1 \in M_1$ and  $n_2 \in N_2$ such that $m_1$ is adjacent to
$n_2$. Let $n_k \in N_k$, then $n_k \in N(m_1)$. By \eqref{path} there exists a
path $P$ from $n_2$ to $n_k$ with $V(P) \subseteq \bigcup_{i=2}^k K_i$ with
$|V(P) \cap K_i|=1$ for every $i \in \{2, \ldots, k\}$.
But now we get a contradiction applying \ref{holenbrs} to the hole $m_1 \dd n_2 \dd P \dd n_k \dd m_1$ and the vertex $v$. This proves~\eqref{claims}.2
 and completes the proof of \eqref{claims}.

  \vspace*{-0.4cm}

\begin{equation} \label{outside}
  \longbox{\emph{Let $v \in V(G) \setminus W$ and 
      for $i \in \{1, \ldots, k\}$, let
  $N_i=N(v) \cap K_i$ and $M_i=K_i \setminus N_i$.
      Either
      \begin{itemize}
      \item $N_i$ is non-empty for at most  two consecutive values of $i$
        (mod $k$) or
      \item $k=5$, $v$ is complete to $W$, and $K_i$ is complete to $K_{i+1}$
        for every $i \in \{1, \ldots, 4\}$.
  \end{itemize}}}
\end{equation}

First we claim that we can choose $j,l$ with $N_j \neq \emptyset$ and $N_l \neq \emptyset$, and such that $j = l \pm 2$.
If $N_i \neq \emptyset$ for every $i \in \{1, \ldots, k\}$, then
the claim holds. If $N_i=\emptyset$ for every $i \in \{1, \ldots, k\}$, then
the claim holds.
Thus we may assume that some $N_i$s are empty, and some are
not. By shifting the indices, we may assume $N_1 \neq \emptyset$ and
$N_k=\emptyset$. We may assume that $N_t \neq \emptyset$ for some
$t \in \{3, \ldots, k-1\}$ for otherwise \eqref{outside} holds with $i=1$.
Let $n_1 \in N_1$ and  $n_t \in N_t$.
By \eqref{path} there exists
a path $P$ from $n_1$ to $n_t$ such that
$V(P) \subseteq \bigcup_{j=1}^t K_j$ and $P$ uses exactly one vertex from each of
$K_1, K_2, \ldots, K_t$.
Also by \eqref{path} there exists a path $Q$ from $n_t$ to $n_1$ 
such that  $V(Q) \subseteq K_1 \cup \bigcup_{j=t}^kK_j$ and $Q$  uses exactly
one vertex from each of $K_t, K_{t+1} \ldots, K_k,K_1$. Now
$F=n_1 \dd P \dd n_t \dd Q \dd n_1$ is a  hole, and
$n_1,n_t \in V(F)$. Since $t \geq 3$ and $N_k=\emptyset$, 
applying \ref{holenbrs} to $F$ and $v$ we deduce
that the fourth outcome of \ref{holenbrs} holds and $t=3$. Now  we can set $j=1$ and $l=t$. This proves the claim.

By the claim of the previous paragraph (shifting the indices so that
$j=k$ and $l=2$) we may
assume that $N_k$ and $N_2$ are both non-empty.
For $i \in \{2,k\}$ let $a_i \in N_i$. By \eqref{complete} $a_k$ and $a_2$ have a common neighbour $a_1 \in K_1$.

    Suppose $\bigcup_{i=3}^{k-1}N_i=\emptyset$. Since $W$ is maximal,
    $(K_1 \cup \{v\}, K_2, \ldots, K_k)$ is not a $k$-structure
    for $W \cup \{v\}$, and therefore there exists $a_1' \in M_1$.
    By symmetry we may assume that $K_1$ is complete to $K_2$.
    Let $a_3' \in K_3$ be a neighbour of $a_2$; then $a_3' \in M_3$,
    contrary to \eqref{claims}.1. This proves that
    $\bigcup_{i=3}^{k-1}N_i \neq \emptyset$.

    Suppose $k=4$. Then there is symmetry between $K_1$ and $K_3$, and
    we deduce that $N_i \neq \emptyset$ for every $i \in \{1, \ldots, 4\}$.
    By \eqref{complete} we may assume that $K_1$ is complete to $K_2$, and
    $K_3$ to $K_4$. Now \ref{holenbrs} implies that there is no hole
    $n_1 \d n_2 \d n_3 \d n_4 \d n_1$ with $n_i \in N_i$ for every
    $i \in \{1, \ldots, 4\}$, and consequently either
    $N_1$ is anticomplete to $N_4$, or $N_3$ is anticomplete to $N_2$.
    By symmetry we may assume that $N_1$ is anticomplete to $N_4$.
    Let $n_1 \in N_1$ and $n_4 \in N_4$, and let $m_1 \in K_1$ be adjacent to
    $n_4$ and $m_4 \in K_4$ be adjacent to $n_1$. Then $m_1 \in M_1$ and
    $m_4 \in M_4$. By \eqref{comparable} applied to $m_1$ and $n_1$, we
    deduce that $m_1$ is adjacent to $m_4$. 
    If there exists $m_2 \in M_2$,
    then $\{n_1,v,m_2,m_4\}$ is a claw, a contradiction. This
    proves (using symmetry) that $M_2 \cup M_3 = \emptyset$. Let $n_2 \in K_2$,
    and let $n_3 \in K_3$ be adjacent to $n_2$. Then $n_2 \in N_2$ and
    $n_3 \in N_3$.  But now $G[v,m_4,n_2,n_4,n_1,n_3,m_1]$ is a seven-antihole,
    a contradiction.  This proves that $k \geq 5$.

By \eqref{claims}.2 it follows that $a_1 \in N_1$.
We claim that $k=5$ and $v$ is complete to $W$.
Suppose $v$ has a non-neighbour in $W$.
Since $\{v,a_k,a_2,x\}$ is not a claw for any
  $x \in \bigcup_{i=4}^{k-2}K_i$, it follows that 
  $\bigcup_{i=4}^{k-2}N_i=\emptyset$.

We may assume that there is  $a_3 \in N_3$. Since
$\{v,a_2,a_k,a_3\}$ is not a claw, it follows that $a_2$ is adjacent to $a_3$.
By \eqref{path} 
  there is a path $P$ from $a_3$ to $a_k$ with $V(P) \subseteq \bigcup_{j=3}^kK_j$
  and using exactly one vertex from each of $K_3, \ldots, K_k$.
  Now $F=a_k \dd a_1 \dd a_2 \dd a_3 \dd P \dd a_k$ is a hole,
  and \ref{holenbrs} implies that $k=5$ and $v$ is complete to $V(F)$.
  We have proved that $N_i \neq \emptyset$ for every $i \in \{1, \ldots, 5\}$
  thus restoring the symmetry of the $5$-structure. Since for  $n_1  \in N_1$,
  $n_2 \in N_2$ and $n_4 \in N_4$, $\{v,n_1,n_2,n_4\}$ is not a claw,
  we deduce (using symmetry) that $N_i$ is complete to $N_{i+1}$ for
  every $i \in \{1, \ldots, 5\}$.

  Next suppose that both $M_5$ and $M_2$ are non-empty. By \eqref{complete} we
  may assume that $K_1$ is complete to $K_2$. By \eqref{claims}.1,
  $N_1$ is anticomplete to $M_5$. Since every
  vertex of $M_5$ has a neighbour in $K_1$, it follows that $M_1 \neq \emptyset$.
By \eqref{claims}.2 
$M_1$ is anticomplete to
  $N_5$,  but now $m_1 \in M_1$ and $n_1 \in N_1$ contradict
  \eqref{comparable}. By symmetry we may assume  $M_i$ is non-empty for
  at most two consecutive values of $i$, and that
  $M_1 \cup M_2 \cup M_3=\emptyset$. 
 By \eqref{claims}.2
 $N_4$ is anticomplete to $M_5$,
   and similarly
   $M_4$ is anticomplete to $N_5$.
By symmetry we may assume that $M_4 \neq \emptyset$.
   But now $m_4 \in M_4$
   and $n_4 \in N_4$ contradicts \eqref{comparable}.
   This proves that $k=5$ and $v$ is complete to $W$.
   To complete the proof of \eqref{outside} assume for a contradiction that
   there exist $i \in \{1, \ldots, 4\}$, $k_i \in K_i$ and $k_{i+1} \in K_{i+1}$
   such that $k_i$ is non-adjacent to $k_{i+1}$. We may assume $i=1$.
   Then $\{v, k_1,k_2,k_4\}$ is a claw in $G$, a contradiction.
   This proves \eqref{outside}.
  \\
  \\
  For $i \in \{1, \ldots, k\}$ let
   $K_{i,i+1}$ be the set of all vertices of $V(G) \setminus W$ that have a neighbour in $K_i$, a neighbour in $K_{i+1}$, and no neighbour in
$W \setminus (K_i \cup K_{i+1})$.
  The  outcomes of \eqref{outside} are summarized in Figure
  \ref{fig:outside}.

  \begin{figure}[h]
\begin{center}
\begin{tikzpicture}[scale=0.28]


\node[inner sep=8pt,   circle, draw] at (32, 4)(v1){}; 
\node[inner sep=8pt,  circle, draw] at (35, 8)(v2){}; 
\node[inner sep=8pt, circle, draw] at (41, 8)(v3){};
\node[inner sep=8pt , circle, draw] at (44, 4)(v4){}; 
\node[inner sep=8pt, circle, draw] at (38, 0)(v5){}; 
\node[inner sep=2.5pt, fill=black, circle] at (38, 4)(v6){};

\node at (32,4) {$K_5$};
\node at (35,8) {$K_1$};
\node at (41,8) {$K_2$};
\node at (44,4) {$K_3$};
\node at (38,0) {$K_4$};
\node at (38.5, 3) {$v$};

\draw[black, thick] (v1) -- (v2);
\draw[black, thick] (v2) -- (v3);
\draw[black, thick] (v3) -- (v4);
\draw[black, thick] (v4) -- (v5);
\draw[black, thick] (v5) -- (v1);

\draw[black, thick] (v1) -- (v6);
\draw[black, thick] (v2) -- (v6);
\draw[black, thick] (v3) -- (v6);
\draw[black, thick] (v4) -- (v6);
\draw[black, thick] (v5) -- (v6);

\node[inner sep=8pt, circle, draw] at (12, 4)(v1){}; 
\node[inner sep=8pt, circle, draw] at (26, 4)(v2){}; 
\node[inner sep=8pt, circle, draw] at (16, 8)(v3){}; 
\node[inner sep=8pt, circle, draw] at (22, 8)(v4){};
\node[inner sep=2.5pt, fill=black, circle] at (19, 4)(v5){}; 

\node at (12,4) {$K_{i-1}$};
\node at (16,8) {$K_{i}$};
\node at (22,8) {$K_{i+1}$};
\node at (26,4) {$K_{i+2}$};
\node at (19,3) {$v$};

\draw[black, dotted, thick] (v1)  .. controls +(0,-7) and +(0,-7) .. (v2);
\draw[black, thick] (v1) -- (v3);
\draw[black, thick] (v2) -- (v4);
\draw[black, thick] (v3) -- (v4);

\tikzset{decoration={snake,amplitude=.4mm,segment length=2mm,
                       post length=0mm,pre length=0mm}}

\draw[decorate] (v3) -- (v5);
\draw[decorate] (v4) -- (v5);

\end{tikzpicture}
\end{center}
\vspace{-0.8cm}
\caption{Outcomes of \eqref{outside} (here wiggly lines represent  possible
  adjacency, and the dotted arc represents the remainder of the $k$-structure).}
\label{fig:outside}
\end{figure}
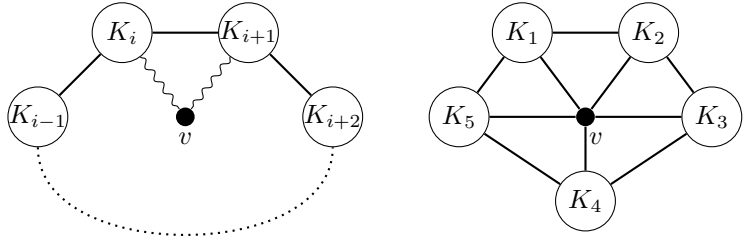

\vspace*{-0.4cm}

\begin{equation} \label{hats}
  \longbox{\emph{Assume that $K_{i,i+} \neq \emptyset$.
      The following statements hold
      \begin{enumerate}
\item      $K_i\cup K_{i+1} \cup K_{i,i+1}$ is a clique.
\item       If $u \in V(G) \setminus W$ is complete to $W$, then $u$ is anticomplete to
  $K_{i,i+1}$.
\item  $K_{i,i+1}$ is anticomplete to $K_{i-1,i}$.
  \end{enumerate}}}
\end{equation}

Let $v \in K_{1,2}$. For $i \in \{1,2\}$ let $N_i=K_i \cap N(v)$, and
let $M_i=K_i \setminus N_i$.  By \eqref{claims}.1 $N_1$ is anticomplete to $M_2$, and $N_2$  is anticomplete to $M_1$. If there exists $m_1 \in M_1$, then
$n_1,m_1$ contradict \eqref{comparable} for every $n_1 \in N_1$.  Thus $M_1=\emptyset$, and by
symmetry $M_2=\emptyset$. This proves that $K_{1,2}$ is complete to
$K_1 \cup K_2$.

Suppose $k_1 \in K_1$ is non-adjacent to $k_2 \in K_2$.
Let $P$ be a path from $k_2$ to a vertex $k_k \in K_k$ as in \eqref{path},
such that $|V(P) \cap K_i|=1$ for every $i \in \{2, 3, \ldots, k\}$.
By \eqref{complete} $K_1$ is complete to $K_k$, and so
$F=k_k \dd k_1 \dd v \dd k_2 \dd P \dd k_k$ is  a hole.
Let $k_1' \in K_1 \cap N(k_2)$, then $G[V(F) \cup \{k_1'\}]$ is a jewel in
$G$, a contradiction. This proves that $K_1$ is complete to $K_2$.

Since $\{k_1,k_k,a,b\}$  is not a claw for any $k_1 \in K_1$,
$k_k \in K_k \cap N(k_1)$ and
$a,b \in K_{1,2}$, it follows that
$K_{1,2}$ is a clique, and therefore $K_1 \cup K_2 \cup K_{1,2}$ is a clique.
By symmetry, $K_1 \cup K_2 \cup K_{1,2}$ is a clique for every $i$.
This proves \eqref{hats}.1.

Next suppose that $u$ is complete to $W$ and  $u$ is adjacent to
$w \in K_{1,2}$.  By \eqref{outside} $k=5$. Let $k_3 \in K_3$ and $k_5 \in K_5$.
Then $\{u,w,k_3,k_5\}$ is a claw, a contradiction. This proves 
\eqref{hats}.2.

Finally, suppose that $w_k \in K_{k,1}$ is adjacent to $w_2 \in K_{1,2}$.
Let $k_1 \in K_1$, $k_2\in K_2$ and $k_k \in K_k$, and let $P$ be a path from $k_2$ to $k_k$ as in \eqref{path} such that $|V(P) \cap K_i|=1$ for every $i \in \{2,3, \ldots, k\}$.
Then $F=k_2 \dd P \dd k_k \dd w_k \dd w_2 \dd k_2$ is a hole
and $G[V(F) \cup \{k_1\}]$ is a jewel, a contradiction. This proves
\eqref{hats}.3 and completes the proof of \eqref{hats}.
\\
\\
For $i \in \{1, \ldots, k\}$, let $V(H) \cap K_i=\{h_i\}$.
Choose $i \in \{1, \ldots, k\}$ such that $K_i$ is complete to
$K_{i+1}$   and, if possible, such that $K_{i,i+1} \neq \emptyset$; we may assume $i=1$. Set $a=h_1$ and $b=h_2$, and let $K=K_1 \cup K_2 \cup K_{1,2}$.

By \eqref{outside}, every  vertex of $X(ab)=\{a,b\} \cup (N(a) \cap N(b))$
either
belongs to $K$ or is complete to $W$ (and $k=5$). Since if
$y \in X(ab)$ is complete to $W$, then $y$ has two non-adjacent
neighbours in $V(G) \setminus X(ab)$, it follows that 
$K$ contains the dome of $ab$.
By \eqref{hats}.1 $K$ is a clique.

We prove that $K$ is a simplicial clique, and therefore $K$ equals the dome
of $ab$.  Suppose $K$ is not a simplicial clique, and let $v \in K$ and
$u,w \in V(G) \setminus K$ be such that
  $u$ and $w$ are adjacent to $v$, and $u$ is non-adjacent to $w$. Suppose first that
$u$ is complete to $W$. By \eqref{outside} $k=5$ and for every $i$
$K_i$ is complete to $K_{i+1}$.
By \eqref{hats}.2 $v \not \in K_{1,2}$.
  For $i \in \{1, \ldots, 5\}$ let $k_i \in K_i$.
  We may assume that $v=k_1$. Since $u$ is non-adjacent to $w$,
  it follows that $w \not \in W$.
  Since $G[k_1,u,k_3,w,k_2]$ is not a jewel in $G$,
    it follows that
    $w$ is not complete to $W$. By \eqref{claims}.1 (since $k_1$ is complete to
    $K_2 \cup K_5$) $w$ has a neighbour in at least one of  $K_2,K_5$. By \eqref{outside} $w \in K_{1,2} \cup K_{5,1}$. Since $w \not \in K$, it follows that $w \in K_{5,1}$.
By \eqref{hats}.1
$K_5$ is complete to $K_1$. Since $K_{1,5} \neq \emptyset$, and since $K_1,K_2$
where chosen with $K_{1,2} \neq \emptyset$ if possible, 
it follows  that
  there exists $k \in K_{1,2}$. By \eqref{hats}.2 $u$ is non-adjacent to $k$,
  and by \eqref{hats}.3 $w$ is non-adjacent to $k$. But now
  $\{k_1,u,k,w\}$ is  a claw in $G$, a contradiction.
    This proves that $u$ is not
  complete to $W$. By symmetry, $w$ is not complete to $W$.

  Now suppose $v \in K_1$. Since $u,w \not \in K$, it
  follows from \eqref{outside} that $u,w \in K_k \cup K_{k,1}$. But then
  $u$ is adjacent to $w$ by \eqref{hats}.1, a contradiction.
  This proves that $v \not \in K_1$, and,
  by symmetry, $v \not \in K_2$.
  
  It follows that $v \in K_{1,2}$. 
Moreover, applying the previous argument to an
  arbitrary vertex of $K_1 \cup K_2$ in the role on $v$, we deduce that
  no vertex of $K_1 \cup K_2$ is complete to $\{u,w\}$.
  Since $\{v,u,w,p\}$ is not a claw for $p \in K_1 \cup K_2$
    it follows that every vertex of  $K_1 \cup K_2$ has a
    neighbour in $\{u,w\}$. Since $u,w \not \in K$, \eqref{hats}.1
    implies that each of $u,w$ is anticomplete to at least one of $K_1,K_2$.
    By switching $u$ and $w$ if necessary, we may assume that
    $u$ has a neighbour in  $K_1$ and is anticomplete to $K_2$.
    By \eqref{claims}.1 $u$ has a neighbour in  $K_k$, but now $u \in K_{k,1}$
    is adjacent to $v \in K_{1,2}$, contrary to \eqref{hats}.3.
Thus we have found an edge of $H$ whose dome is a simplicial clique.
    This proves \ref{subgraphs}.
    \end{proof}

\section{The Algorithm}

In this section we use \ref{subgraphs} to design a simple algorithm
that finds a simplicial clique in a clean graph.

\begin{theorem} \label{alg}
  There is an algorithm with the following specifications.\\
  {\bf Input:} A non-null clean graph $G$.\\
  {\bf Output:} A simplicial clique in $G$.\\
  {\bf Running time:} $O(|V(G)|^5)$.
\end{theorem}

\begin{proof}
  Let $|V(G)|=n$.
  First, for every vertex $v \in V(G)$ check if $N(v)$ is a clique.
  If the answer is yes for some $v$, then output a simplicial clique $\{v\}$.
  This step can be done in time $O(n^3)$.

  Now for every edge $ab$ compute $X(ab), Y(ab)$ and the dome of $ab$,
  and check if the dome of $ab$ is a simplicial clique.
     This step can be done in time $O(n^5)$.

     We now use \ref{subgraphs} to  prove that the algorithm will return a
     simplicial clique of $G$. If $G$ is a chordal graph, then
     by the first statement of \ref{subgraphs} the first step of the algorithm
     will return a simplicial clique; otherwise there is a hole in $G$, and so
     by the second statement of \ref{subgraphs}, the second step of the algorithm will return a simplicial clique.
     This proves \ref{alg}.
\end{proof}

We remark that the algorithm of \ref{alg} can be used to produce, in
time $O(|V(G)|^2)$, a list of at most $|V(G)|^2$ sets one of which is guaranteed to be a simplicial clique. The rest of the running time is spent checking whether each of the sets is a simplicial clique.

\section{Acknowledgment}

We are grateful to Adrian Chapman and Steve Flammia for their help in writing
the introduction.

\end{document}